\pdfoutput=1

\documentclass[
    10pt,
    final,
    journal,
    letterpaper,
    twoside,
    twocolumn
]{IEEEtran}
\usepackage[pass]{geometry}
\usepackage{amsmath, amssymb}
\usepackage{array}
\usepackage{url}
\usepackage{booktabs}
\usepackage{siunitx}
\usepackage{cite}
\usepackage{graphicx}
\usepackage[pdfborder={0 0 0}, hidelinks]{hyperref}
\usepackage[capitalize]{cleveref}
\usepackage{textcomp}
\usepackage{microtype}

\DeclareGraphicsExtensions{.pdf}
\graphicspath{{./pdf/}}

\crefformat{equation}{(#2#1#3)}
\crefname{figure}{Fig.\!}{Fig.\!}
\Crefname{figure}{Fig.\!}{Fig.\!}

\newcommand*{\ie}{\textit{i}.\textit{e}.}  
\newcommand*{\eg}{\textit{e}.\textit{g}.}  
\newcommand*{\etal}{\textit{et al}.}       

\newcommand*{\minus}[2][\empty]{{#2}^{{#1}\scriptscriptstyle{-}}}
\newcommand*{\plus}[2][\empty]{{#2}^{{#1}\scriptscriptstyle{+}}}
\newcommand*{\uvar}[1]{\vphantom{#1}\smash{\underline #1}}
\newcommand*{\ovar}[1]{\vphantom{#1}\smash{\overline #1}}
\newcommand*{\<}{\mkern-0mu}

\let\originalleft\left
\let\originalright\right
\renewcommand{\left}{\mathopen{}\mathclose\bgroup\originalleft}
\renewcommand{\right}{\aftergroup\egroup\originalright}

\AtBeginDocument{\DeclareSIUnit{\kWh}{kWh}}
\AtBeginDocument{\DeclareSIUnit{\MWh}{MWh}}
\AtBeginDocument{\DeclareSIUnit{\pu}{pu}}
\begin{document}
\title{Sharing Energy Storage Between\protect\\
       Transmission and Distribution}
\author{Ryan~T.~Elliott,
        Ricardo~Fern\'{a}ndez-Blanco,
        Kelly~Kozdras,
        Josh~Kaplan,
        Brian~Lockyear,
        Jason Zyskowski,
        and~Daniel~S.~Kirschen%
\thanks{R. T. Elliott, K. Kozdras, and D. S. Kirschen
are with the University of Washington Department of Electrical
Engineering in Seattle, WA 98195 USA (e-mail: ryanelliott@ieee.org).}%
\thanks{R. Fern\'{a}ndez-Blanco is with the Department of Applied
Mathematics at the University of Malaga, Malaga 29076 Spain
(e-mail: ricardo.fcarramolino@uma.es).}%
\thanks{J. Kaplan and B. Lockyear are with Doosan GridTech in Seattle, WA.}%
\thanks{J. Zyskowski is with Snohomish County PUD in Everett, WA.}%
}

\IEEEaftertitletext{\vspace{-1\baselineskip}}

%



\maketitle

\begin{abstract}
This paper addresses the problem of how best to coordinate, or
``stack,'' energy storage services in systems that lack centralized
markets.  Specifically, its focus is on how to coordinate
transmission-level congestion relief with local, distribution-level
objectives.  We describe and demonstrate a unified communication and
optimization framework for performing this coordination.  The
congestion relief problem formulation employs a weighted
${\ell_{1}\text{-norm}}$ objective.  This approach determines a set of
corrective actions, \ie, energy storage injections and conventional
generation adjustments, that minimize the required deviations from a
planned schedule.  To exercise this coordination framework, we present
two case studies.  The first is based on a 3-bus test system, and the
second on a realistic representation of the Pacific Northwest region
of the United States. The results indicate that the scheduling
methodology provides congestion relief, cost savings, and improved
renewable energy integration.  The large-scale case study informed the
design of a live demonstration carried out in partnership with the
University of Washington, Doosan GridTech, Snohomish County PUD, and
the Bonneville Power Administration.  The goal of the demonstration
was to test the feasibility of the scheduling framework in a
production environment with real-world energy storage assets.  The
demonstration results were consistent with computational simulations.
\end{abstract}


%
\IEEEpeerreviewmaketitle

\section{Introduction}
\IEEEPARstart{U}{tility-scale} energy storage has the potential to
provide non-wire solutions to longstanding power grid problems.
For example, distribution system operators (DSOs) could use energy
storage to help reduce energy imbalance expenses or to serve their
load more economically through energy arbitrage.
Likewise, transmission system operators (TSOs) could use energy
storage to mitigate congestion or provide frequency regulation.
While the prospect of employing energy storage to tackle these
challenges has drawn immense interest, the question of how best to
coordinate, or ``stack,'' services remains open.
A systematic approach to service coordination would not only help
maximize resource utilization, but also bolster the financial
viability of energy storage projects.

In this work, we describe and demonstrate a unified communication and
optimization framework for scheduling multiple simultaneous storage
services between a TSO and one or more DSOs.  To address this problem,
we propose a multistage approach based on mixed-integer linear
programming.
To demonstrate the viability of the framework, it was implemented and
used to control a utility-scale battery energy storage system (ESS) in
Everett, WA.  This battery is owned and operated by Snohomish County
PUD (SnoPUD), and offers services to the Bonneville Power
Administration (BPA).

The scheduling framework presented here reflects the operating
environment of the Pacific Northwest region of the United States;
however, it is also suitable for systems that lack centralized markets
or that rely heavily upon bilateral contracts. This is the case, for
example, in large parts of Europe~\cite{Hob:04}.  Furthermore, the
ideas developed in this paper may inspire new approaches to energy
storage scheduling in systems that do have centralized markets.

\subsection{Literature review}
\label{sub:lit_review}

The concept of providing multiple simultaneous services with storage
resources has generated active discussion throughout academia,
industry, and government.  The business case for multiple service
provision in wholesale and retail markets is assessed by Teng and
Strbac in \cite{Ten:16}.  Related economic analyses can be found in
\cite{He:11, Hu:16, Fon:17}.  In contrast to~\cite{Ten:16}, which
focuses on the aggregation of distributed storage systems, we
concentrate on independently scheduled utility-scale systems.

The problem of scheduling multiple services is addressed in
\cite{Meg:15, You:10, Gan:14, Kim:17}.  For previous work on sharing
storage resources among multiple parties, see \cite{Par:15, Tus:16,
Yao:16, Kal:17}.  In \cite{Meg:15}, M\'{e}gel \etal\ employ a model
predictive control approach for co-optimizing simultaneous provision
of local and system-wide services.  The algorithm proposed therein
determines the optimal power and energy capacity to allocate for each
service as a function of time.  Alternatively, we take a decentralized
approach that permits the resource to determine the capacity required
for local service provision.  In regard to economics, \cite{Meg:15}
shows that stacking services can improve the financial prospects of
storage resources.

Providing transmission-level congestion relief with energy storage is
explored in \cite{Hu:11, Del:14, Var:15, Kha:16}.  The related problem
of employing energy storage in congestion-constrained distribution
networks is considered in \cite{Koe:04}.  The multi-objective
formulation developed by Khani \etal\ in \cite{Kha:16} seeks to
maximize ESS revenue generated via arbitrage and the ESS contribution
to congestion relief.  Balancing the two objectives is achieved using
an adaptive penalty mechanism.  This mechanism has a similar
mathematical structure to the weighted ${\ell_{1}\text{-norm}}$
employed in the formulation we present in \cref{sec:formulation}.

Regulatory agencies and independent system operators have taken steps
to facilitate the integration of storage resources into markets for
electricity and ancillary services.  In a notice of proposed
rulemaking~\cite{FERC:16}, FERC states that permitting storage
resources to manage their own state of charge would
``allow these resources to optimize their operations to provide
all of the services that they are technically capable of providing.''
We have adopted this perspective in the service coordination framework
developed in this paper.  In line with FERC's guidance, CAISO has held
workshops on the governance of storage resources for multiple-use
applications \cite{CAISO:17}.

\subsection{Paper organization}
\label{sub:paper_organization}

The remainder of this paper is organized as follows.
\cref{sec:methodology} describes the method behind the linked
communication and optimization procedures.  It also outlines the
reports that facilitate communication between parties.  The
formulation of the congestion relief optimization problem is presented
in \cref{sec:formulation}.  In Sections \ref{sec:3bus_case} and
\ref{sec:case_study_and_demo}, we discuss results from the case
studies and live demonstration.  \cref{sec:conclusion} summarizes and
concludes.

\section{Proposed method}
\label{sec:methodology}

\begin{table}[!b]
    \renewcommand{\arraystretch}{1.15}
    \caption{Description of the reports}
    \label{tab:report_description}
    \centering
    \begin{tabular}{l l}
        \toprule
        Name & Contents \\
        \midrule
        Capacity & ESS power and energy capacity \\
        Congestion forecast & Load forecast and TEPO charging indicator \\
        Initial schedule & Initial DEPO injection schedule \\
        Mitigation needs & Minimum and maximum net load \\
        Final schedule & Final ESS combined injection schedule \\
        \bottomrule
    \end{tabular}
\end{table}


The framework developed in this paper is an interlinked series of
optimization problems and data transfers.  We refer to the
transmission grid energy positioning optimizer as TEPO and its
counterpart in the distribution grid as DEPO.  TEPO utilizes available
energy storage capacity to satisfy transmission-side objectives, and
DEPO seeks to satisfy local, distribution-side objectives.

The communication between TEPO and a given DEPO instance is based on
five \textit{reports} that are exchanged when their contents are
required by a particular stage of the optimization.
\cref{tab:report_description} outlines the contents of these reports,
and \cref{tab:info_exchange} the sequence in which they are exchanged.
This scheme complies with the OpenADR specification, an open and
interoperable information exchange model for smart grid
applications~\cite{Sam:16, McP:11}.  In its simplest form, the
procedure in \cref{tab:info_exchange} represents an exchange between
two parties; however, there is no restriction on the number of DEPO
instances with which TEPO may communicate.  Each DEPO instance
corresponds to a distribution-side entity or a subset of a given
entity's storage resources.  In this way, the framework accommodates
service coordination between a transmission system operator and an
arbitrary number of additional parties.



\begin{table}[!t]
    \renewcommand{\arraystretch}{1.15}
    \caption{Communication procedure}
    \label{tab:info_exchange}
    \centering
    \begin{tabular}{c l}
        \toprule
        Direction & Capacity exchange \\
        \midrule
        TEPO$\,\rightarrow\,$DEPO & TEPO requests the capacity report. \\
        TEPO$\,\leftarrow\,$DEPO & DEPO returns the capacity report. \\[1ex]
        \midrule
        Direction & Congestion forecast exchange \\
        \midrule
        TEPO$\,\rightarrow\,$DEPO & TEPO requests the initial schedule report. \\
        TEPO$\,\leftarrow\,$DEPO & DEPO requests the congestion forecast report. \\
        TEPO$\,\rightarrow\,$DEPO & TEPO returns the congestion forecast report. \\
        TEPO$\,\leftarrow\,$DEPO & DEPO returns the initial schedule report. \\[1ex]
        \midrule
        Direction & Mitigation needs exchange \\
        \midrule
        TEPO$\,\rightarrow\,$DEPO & TEPO requests the final schedule report. \\
        TEPO$\,\leftarrow\,$DEPO & DEPO requests the mitigation needs report. \\
        TEPO$\,\rightarrow\,$DEPO & TEPO returns the mitigation needs report. \\
        TEPO$\,\leftarrow\,$DEPO & DEPO returns the final schedule report. \\
        \bottomrule
    \end{tabular}
\end{table}

\begin{figure}[!t]
    \centering
    \includegraphics[width=2.6in]{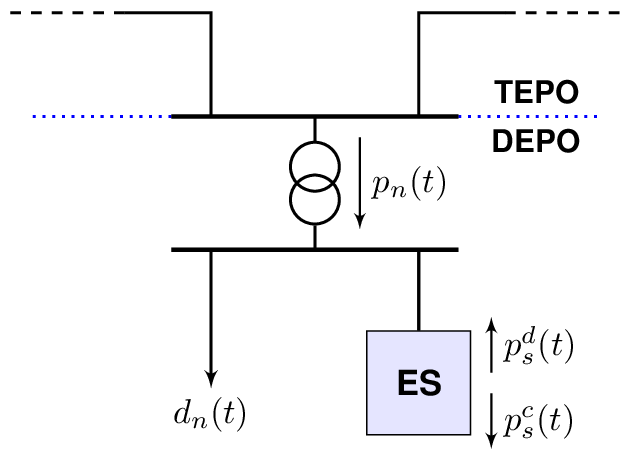}
    \caption{Representative interconnection diagram for an ESS.}
    \label{fig:bus_diagram}
\end{figure}

Here we explain the chain of data transfers for the basic case with
one DEPO instance corresponding to a single ESS.
\Cref{fig:bus_diagram} shows the interconnection of a typical ESS with
the transmission grid.  \Cref{fig:day_ahead_flow} shows the
relationship between the communication procedure and the TEPO
formulation.  Initiating the procedure, TEPO requests the
\textit{Capacity report} from DEPO.  This prompts DEPO to share
information about the power rating and energy capacity of the ESS.
Based on this information and the anticipated system operating
conditions, TEPO provides DEPO with the \textit{Congestion forecast
report}.  This report contains a forecast of the demand at the bus
where the ESS is located and a charging indicator that flags whether
TEPO would like to charge or discharge the ESS in each period of the
optimization horizon.  DEPO uses this information to generate a
preliminary ESS schedule that it shares in the \textit{Initial
schedule report}.  TEPO then processes DEPO's initial schedule and
generates its preferred supplemental injections.  These are
transmitted to DEPO in the form of bounds on the net load at the
energy storage bus in the \textit{Mitigation needs report}.  After
receiving the net load bounds, DEPO finalizes the ESS schedule and
notifies TEPO through the \textit{Final schedule report}.  This
concludes the procedure.


\subsection{Hour-ahead framework}
The scheduling method presented here comprises both day-ahead and
hour-ahead components.  The purpose of the hour-ahead framework is to
reassess the day-ahead schedule and account for changes in the system
operating conditions, \eg, the load and renewable energy forecasts.
This reassessment reflects the fact that there is less uncertainty in
the hour-ahead framework than in the day-ahead.  The flow of
information is effectively the same as in \cref{fig:day_ahead_flow},
except the routine is solved iteratively over a shorter time horizon.
%

\section{Formulation}
\label{sec:formulation}

\begin{figure}[!t]
    \centering
    \includegraphics[width=3.25in]{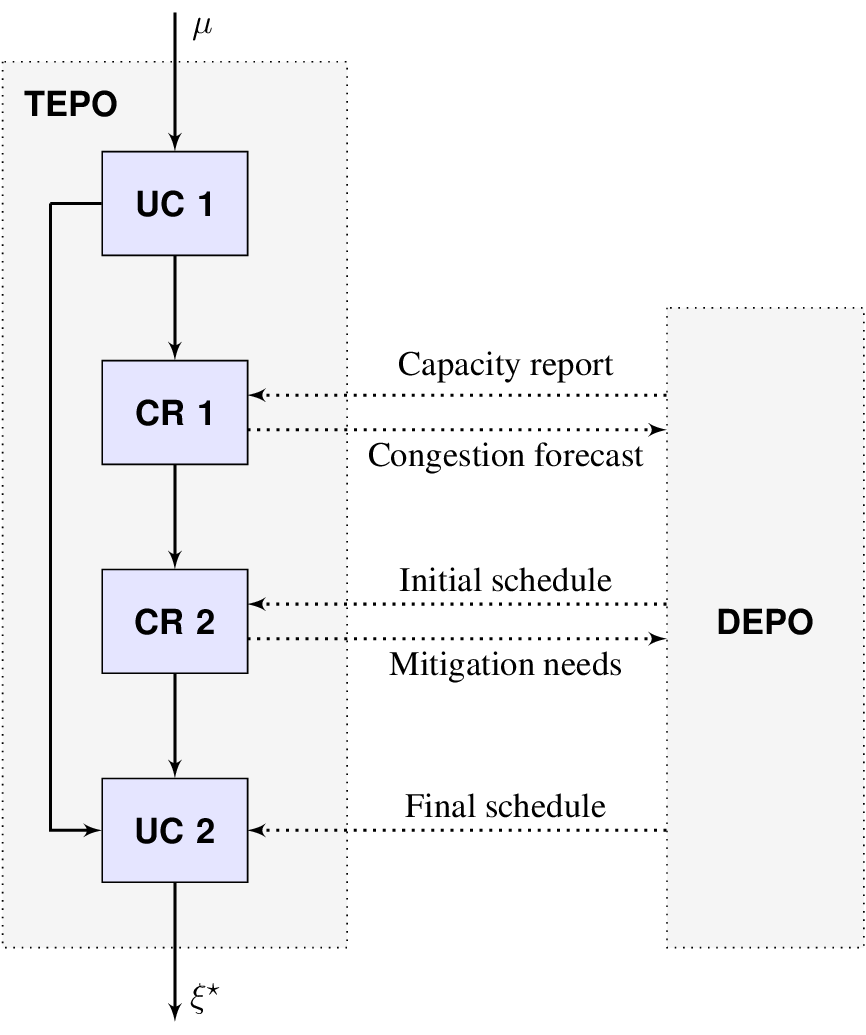}
    \caption{Day-ahead communication and optimization flow.}
    \label{fig:day_ahead_flow}
\end{figure}

The exchange of information delineated in \cref{tab:info_exchange} is
designed to support a multistage optimization framework.  Here we
envision a TEPO formulation that provides transmission-level
congestion relief; however, other services, such as frequency
regulation could also be handled under this framework.  DEPO has the
flexibility to run a separate optimization or control scheme that
oversees local service provision.  \cref{tab:set_nomenclature}
outlines a set nomenclature for the TEPO formulation, and
\cref{tab:stage_description} the optimization stages.  The
relationship between the TEPO formulation and the communication
procedure is shown in \cref{fig:day_ahead_flow}.  At a high level,
each stage of the formulation can be stated as
\begin{equation}
    \label{eq:generic_formulation}
    \begin{array}{l l l}
        \displaystyle{\mathop{\mathrm{minimize}}_{\xi}} & f(\xi), \\
        \mathrm{subject\ to} & \mkern+2mu g_{r}(\xi) \leq 0, & r \in \left\{1, \dots, p\right\}, \\
        & h_{r}(\xi) = 0, & r \in \left\{1, \dots, q\right\},
    \end{array}
\end{equation}
where $\xi$ denotes an ordered list, or \textit{tuple}, of decision
variables.  The inequality constraints are denoted by
${\{g_{r}\}_{r=1}^{p}}$, and the equality constraints by
${\{h_{r}\}_{r=1}^{q}}$.


\subsection{Stage~\num{1}, Pre-mitigation unit commitment}
\label{sec:stage_1}

\begin{table}[!t]
    \renewcommand{\arraystretch}{1.15}
    \caption{Set nomenclature}
    \label{tab:set_nomenclature}
    \centering
    \begin{tabular}{l l l}
        \toprule
        Set & Index & Description \\
        \midrule
        $\mathcal{B}$ & $b$ & Generator cost curve segments \\
        $\mathcal{I}$ & $i$ & Conventional generators \\
        $\mathcal{J}$ & $j$ & Fixed generators \\
        $\mathcal{K}$ & $k$ & Solar power plants \\
        $\mathcal{L}$ & $l$ & Transmission lines \\
        $\mathcal{N}$ & $n$ & Buses \\
        $\mathcal{S}$ & $s$ & Storage devices \\
        $\mathcal{T}$ & $t$ & Time intervals \\
        $\mathcal{W}$ & $w$ & Wind farms \\
        \bottomrule
    \end{tabular}
\end{table}

Stage~\num{1} is called the \textit{Pre-mitigation unit
commitment} problem.
In it, TEPO solves a standard unit commitment and economic
dispatch without taking energy storage or system security
constraints into account.
The output corresponds to the optimal commitment and dispatch
irrespective of energy storage and transmission capacity.

\subsubsection{Decision variables}
The tuple of decision variables in Stage~\num{1} is
\begin{equation}
    \nonumber
    \xi : \left(v_{i}, y_{i}, z_{i}, p_{ib},
        p_{i}, x_{k}, x_{w}, \theta_{n} \right),
\end{equation}
for all ${b \in \mathcal{B}}$, ${i \in \mathcal{I}}$, ${k \in
\mathcal{K}}$, ${n \in \mathcal{N}}$, ${w \in \mathcal{W}}$, and ${t
\in \mathcal{T}}$, where the relevant sets are defined in
\cref{tab:set_nomenclature}.  For the $i$th conventional generation
unit, $v_{i}$, $y_{i}$, and $z_{i}$ are the commitment, startup, and
shutdown status variables, respectively.  Correspondingly, $p_{i}$ is
the total power output, and $p_{ib}$ is the power of the $b$th
segment, or block, of its cost curve.  For renewable generation,
$x_{k}$ is the power curtailment of the $k$th solar plant, and $x_{w}$
the curtailment of the $w$th wind farm.  Lastly, $\theta_{n}$ denotes
the voltage angle at bus~$n$.



\begin{table}[!t]
    \renewcommand{\arraystretch}{1.15}
    \caption{Optimization stages}
    \label{tab:stage_description}
    \centering
    \begin{tabular}{l l l}
        \toprule
        No. & Stage & Description \\
        \midrule
        1 & \textsf{\scriptsize UC 1} & Pre-mitigation stage \\
        2 & \textsf{\scriptsize CR 1} & Independent congestion relief \\
        3 & \textsf{\scriptsize CR 2} & Coordinated congestion relief \\
        4 & \textsf{\scriptsize UC 2} & Post-mitigation stage \\
        \bottomrule
    \end{tabular}
\end{table}


\subsubsection{Objective function}

For Stage~\num{1}, we employ a standard minimum operating cost objective
that can be separated into three components:
\begin{align}
    \label{eq:stage1_objective}
    f(\xi) = \sum_{t \in \mathcal{T}}\sum_{i \in \mathcal{I}}C_{i}(t)
    + \sum_{t \in \mathcal{T}}\sum_{w \in \mathcal{W}}C_{w}(t)
    + \sum_{t \in \mathcal{T}}\sum_{k \in \mathcal{K}}C_{k}(t),
\end{align}
where $C_{i}$ is the total cost incurred by the $i$th conventional
generating unit, while $C_{w}$ and $C_{k}$ account for the cost of
curtailing wind and solar generation, respectively.

The conventional generation costs are given by
\begin{align}
\label{eq:obj_conv_stage1}
C_{i}(t) = c_{i}^{\mathrm{nl}} v_{i}^{}(t)
    + c_{i}^{\mathrm{su}} y_{i}^{}(t)
    + \sum_{b\in\mathcal{B}}m_{ib} \mkern+1mu p_{ib}(t),
\end{align}
where $m_{ib}$ is the incremental cost
of the $b$th block of unit~$i$'s cost curve,
$c_{i}^{\mathrm{nl}}$ is the no-load cost, and
$c_{i}^{\mathrm{su}}$ the startup cost.

The costs arising from curtailing renewable generation are
\begin{align}
    \label{eq:wind_curtail}
    C_{w}(t) &= m_{w} x_{w}(t) \\
    \label{eq:solar_curtail}
    C_{k}(t) &= m_{k} x_{k}(t),
\end{align}
where $m_{w}$ and $m_{k}$ denote the incremental costs of curtailing
the $w$th wind farm and $k$th solar plant, respectively.
%
%

\subsubsection{Binary variable generation constraints}

The startup and shutdown of unit~$i$ is captured by
the constraints
\begin{align}
    \label{eq:binary_one}
    & y_{i}(t) - z_{i}(t) = v_{i}(t) - v_{i}(t - 1) \\
    \label{eq:binary_two}
    & y_{i}(t) + z_{i}(t) \leq 1,
\end{align}
which are enforced for all ${i \in \mathcal{I}}$ and ${t \in \mathcal{T}}$.
In the initial time period, ${v_{i}(t - 1)}$ takes a special
value, $v_{i}^{0}$, that reflects the initial commitment
status of unit~$i$.

\subsubsection{Minimum up and down time constraints}

Let $\ovar{L}_{i}$ denote the number of time periods
that unit~$i$ must remain up at the beginning of
the operating horizon, and let $\uvar{L}_{i}$ be the
corresponding number of periods that it must remain down.
These parameters are defined in \cite{Pan:15}.
At least one of $\ovar{L}_{i}$ and $\uvar{L}_{i}$ is zero
by definition.
%
%
%
To ensure that the commitment status of unit~$i$
remains unchanged for the initial number of periods
dictated by $\ovar{L}_{i}$ or $\uvar{L}_{i}$,
the constraint
\begin{equation}
\label{eq:min_updn_one}
v_{i}(t) = v_{i}^{0}, \text{ for all } t \leq
        \ovar{L}_{i} + \uvar{L}_{i}, \\
\end{equation}
is enforced for all ${i \in \mathcal{I}}$.
The minimum up and down time requirements
over the remainder of the operating horizon are given by
\begin{align}
    \label{eq:min_updn_two}
    \sum_{\tau = t - \ovar{\Gamma}_{i} + 1}^{t} & y_{i}(\tau) \leq v_{i}(t),
        \text{ for all } t \geq \ovar{L}_{i} \\
    \label{eq:min_updn_three}
    \sum_{\tau = t - \uvar{\Gamma}_{i} + 1}^{t} & z_{i}(\tau) \leq 1 - v_{i}(t),
        \text{ for all } t \geq \uvar{L}_{i},
\end{align}
for all ${i \in \mathcal{I}}$, where
$\ovar{\Gamma}_{i}$ is the minimum up time,
and $\uvar{\Gamma}_{i}$ the minimum down time.

\subsubsection{Generator output constraints}

The total power output of the $i$th unit is given by sum of the
outputs corresponding to the segments, or blocks, of its cost curve,
\ie,
\begin{align}
    \label{eq:gen_con_1}
    p_{i}(t) = \sum_{b\in \mathcal{B}} p_{ib}(t),
\end{align}
for all ${i \in \mathcal{I}}$ and ${t \in \mathcal{T}}$.
The conventional generator cost curves employed in
\cref{eq:obj_conv_stage1} are piecewise-linear with segments
${b\in\mathcal{B}}$.  Typically, the associated marginal cost curves
are monotonically nondecreasing.

The power output of each block and the total output of each unit
is bounded such that
\begin{align}
    \label{eq:gen_con_2}
    \uvar{p}_{i}v_{i}(t) &\leq p_{i}(t) \leq \ovar{p}_{i}v_{i}(t) \\
    \label{eq:gen_con_3}
    0 &\leq p_{ib}(t) \leq \ovar{p}_{ib}v_{i}(t),
\end{align}
for all ${b \in \mathcal{B}}$, ${i \in \mathcal{I}}$,
and ${t \in \mathcal{T}}$.
For the $i$th unit, $\uvar{p}_{i}$ is the minimum power output,
and $\ovar{p}_{i}$ the maximum power output.
Similarly, $\ovar{p}_{ib}$ is the maximum output of the
$b$th block of unit~$i$'s cost curve.
The maximum unit output, $\ovar{p}_{i}$, is permitted
to vary between operating horizons to account for
scheduled maintenance and forced outages.

\subsubsection{Ramping constraints}

The ramping constraints on conventional generation
are expressed as follows:
\begin{align}
    \label{st1_ramp_constraints}
    -\uvar{R}_{i} \leq p_{i}(t) - p_{i}(t-1) \leq \ovar{R}_{i},
\end{align}
for all ${i\in\mathcal{I}}$ and ${t\in\mathcal{T}}$,
where $\ovar{R}_{i}$ is the upward ramp limit,
and $\uvar{R}_{i}$ the downward ramp limit.
In the initial time period, $p_{i}(t-1)$ takes a special
value, $p_{i}^{0}$, that reflects the initial power output of
the $i$th unit.

\subsubsection{Renewable generation curtailment constraints}
Recall that the objective function stated in \cref{eq:stage1_objective}
includes terms that account for the cost of curtailing variable
renewable generation.
In accordance with \cref{eq:wind_curtail} and \cref{eq:solar_curtail},
the constraints on curtailment are given by
\begin{align}
    \label{eq:st1_wind_curtail_con}
    0 &\leq x_{w}(t) \leq p_{w}(t) \\
    \label{eq:st1_solar_curtail_con}
    0 &\leq x_{k}(t) \leq p_{k}(t),
\end{align}
where $p_{w}$ is the power available at the $w$th wind farm, and
$p_{k}$ the power available at the $k$th solar plant.  These
constraints are enforced for all ${w \in \mathcal{W}}$,
${k \in \mathcal{K}}$, and ${t \in \mathcal{T}}$.

\subsubsection{Nodal power balance constraints}

Power transfer throughout the transmission network
is modeled using a standard dc power flow approximation, \ie,
\begin{align}
    \label{eq:dc_line_flow}
    F_{l}(t) = \left\lbrack\theta_{o(l)}\<(t) - \theta_{d(l)}\<(t)\right\rbrack/x_{l},
\end{align}
for all ${l \in \mathcal{L}}$ and ${t \in \mathcal{T}}$.  This network
model encompasses high-voltage transmission lines for which the
resistance to reactance ratio, $r_{l}/x_{l}$, may be reasonably
assumed to be small~\cite{Sto:09,Con:02,Mot:02}.  For line $l$,
$x_{l}$ is the reactance, and $F_{l}$ the real power flow.  The
function $o(l)$ returns the origin or ``from'' bus index of line $l$,
and $d(l)$ returns the destination or ``to'' bus index.  The voltage
angles are bounded between $-\pi$ and $\pi$ for all ${n \in
\mathcal{N}}$ and ${t \in \mathcal{T}}$.  Per convention, the
reference bus is constrained to have a voltage angle of zero for all
${t \in \mathcal{T}}$.

Using \cref{eq:dc_line_flow}, the nodal power balance constraints
can be stated as follows:
\begin{multline}
    \label{eq:nodal_power_bal}
    d_{n}(t) = \sum_{i\in\mathcal{I}_{n}} p_{i}(t)
         + \sum_{j\in\mathcal{J}_{n}} p_{j}(t)
         + \sum_{k\in\mathcal{K}_{n}} \hat{p}_{k}(t) \\
         + \sum_{w\in\mathcal{W}_{n}} \hat{p}_{w}(t)
         - \sum_{l\in\mathcal{O}_{n}}F_{l}(t) + \sum_{l\in\mathcal{D}_{n}}F_{l}(t),
\end{multline}
for all ${n \in \mathcal{N}}$ and ${t \in \mathcal{T}}$,
where $d_{n}$ is the total demand at bus~$n$.
Let a subscript $n$ affixed to a set indicate the subset of
components connected to bus~$n$.
The set of lines originating at bus~$n$ is denoted by
$\mathcal{O}_{n}$, and the set with destinations at bus~$n$
by $\mathcal{D}_{n}$.
%
%
%
The net power from the $k$th
solar plant is ${\hat{p}_{k}(t) = p_{k}(t) - x_{k}(t)}$,
and the net power from the $w$th wind farm
${\hat{p}_{w}(t) = p_{w}(t) - x_{w}(t)}$.


\subsection{Stage~\num{2}, Independent congestion relief}

Stage~\num{2} is called the
\textit{Independent congestion relief} problem
because TEPO solves it without knowledge of the ESS schedule
(or hour-ahead adjustments) that DEPO would like to carry out.
Using the capacity report, TEPO solves an optimization to determine
a minimal set of corrective actions, \ie, commitment and dispatch
adjustments and ESS injections, required to alleviate congestion.
This allows TEPO to form a cursory schedule
indicating whether each ESS should be charged or discharged
as a function of time to mitigate congestion.

At a high level, the objective of the congestion relief
problem takes the form
\begin{align}
    \label{eq:weighted_l1}
    f(\xi) = \sum_{r\in\mathcal{R}}\lvert\alpha_{r}\xi_{r}\rvert =
        {\left\lVert A \xi \right\rVert}_{1},
\end{align}
where ${\alpha_{r} \ge 0}$ for all ${r \in \mathcal{R}}$.
Let ${A \in \mathbb{R}^{\lvert\mathcal{R}\rvert \times \lvert\mathcal{R}\rvert}}$
be a real, positive-semidefinite diagonal matrix.
Mathematically, \cref{eq:weighted_l1} is a weighted ${\ell_{1}\text{-norm}}$
where $\alpha_{r}$ is the weight for the $r$th entry of $\xi$.
Recall that the \textit{cardinality} of $\xi$ is the
number of nonzero entries it contains.
Although cardinality minimization is NP-hard in general,
for bounded systems of linear equalities and inequalities
it is equivalent to ${\ell_{p}\text{-norm}}$ minimization \cite{Fun:11}.
As shown in \cite{Don:06}, the ${\ell_{1}\text{-norm}}$
is the convex envelope, \ie, the best convex lower bound,
of the cardinality function.
For this reason, the ${\ell_{1}\text{-norm}}$ is used as a convex
approximation to the cardinality function in statistical
regression, compressed sensing, and elsewhere
\cite{Don:08, Fig:07}.
Hence, the objective stated in (\ref{eq:weighted_l1}) has the
effect of minimizing the \textit{number} of corrective actions
required to alleviate congestion.

It is possible to formulate optimization problems with
${\ell_{1}\text{-norm}}$ objectives as linear programs,
as described in \cite{Man:14}.
Absolute value terms, such as
${\left\lvert\alpha_{r}\xi_{r}\right\rvert}$,
can be implemented with auxiliary variables
and constraints of the form
\begin{align}
    \label{eq:abs_val_1}
    \xi_{r} &= \plus{\xi}_{r} - \minus{\xi}_{r} \\
    \label{eq:abs_val_2}
    0 &\leq \plus{\xi}_{r}, \ \ 0 \leq \minus{\xi}_{r}.
\end{align}
The absolute value is then given by
\begin{equation}
    \label{eq:abs_val_3}
    \left\lvert\alpha_{r}\xi_{r}\right\rvert =
        \lvert\alpha_{r}\rvert\left(\plus{\xi}_{r} + \minus{\xi}_{r}\right),
\end{equation}
or simply ${\alpha_{r}(\plus{\xi}_{r} + \minus{\xi}_{r})}$
where ${\alpha_{r} \ge 0}$.
Where necessary, binary variable constraints can be introduced
to ensure that at most one of $\plus{\xi}_{r}$ and $\minus{\xi}_{r}$
is nonzero.
%
%

\subsubsection{Decision variables}
The tuple of decision variables in Stage~\num{2} is
\begin{align}
    \nonumber
    \xi : \bigl(v_{i}, y_{i}, z_{i}, \nu_{s}, \psi_{s}, \zeta_{s}, E_{s},
        \delta p_{ib},
        \delta p_{i},
        \delta p_{s}^{c},
        \delta p_{s}^{d},
        \delta x_{w},
        \delta x_{k} \bigr),
\end{align}
for all ${b \in \mathcal{B}}$, ${i \in \mathcal{I}}$,
${k \in \mathcal{K}}$, ${n \in \mathcal{N}}$, ${s \in \mathcal{S}}$,
${w \in \mathcal{W}}$, and ${t \in \mathcal{T}}$.
The conventional generation binary variables
are defined as in Stage~\num{1}.
The energy storage binary variables $\nu_{s}$, $\psi_{s}$, and
$\zeta_{s}$ prevent simultaneous charging and discharging
(and simultaneous opposing adjustments).
The state of charge (SOC) of the $s$th ESS is denoted by $E_{s}$.
Let variables beginning with $\delta$ denote
adjustments or deviations in some underlying quantity.
For example, $\delta p_{s}^{c}$ is the charging adjustment
of the $s$th ESS, and $\delta p_{s}^{d}$ the
corresponding discharging adjustment.
Each of these deviations is implemented with
a pair of nonnegative decision variables as in
\cref{eq:abs_val_1}--\cref{eq:abs_val_3}.%

\subsubsection{Objective function}
In this case, the objective function does not correspond
precisely to an economic cost.
Rather, it represents a penalty
for deviating from the schedule determined in Stage~\num{1}.
For Stage~\num{2}, the objective function can be separated into four
components:
\begin{multline}
    \label{eq:stage2_objective}
    f(\xi) = \sum_{t \in \mathcal{T}}\sum_{i \in \mathcal{I}}\Phi_{i}(t)
    + \sum_{t \in \mathcal{T}}\sum_{w \in \mathcal{W}}\Phi_{w}(t) \\
    + \sum_{t \in \mathcal{T}}\sum_{k \in \mathcal{K}}\Phi_{k}(t)
    + \sum_{t \in \mathcal{T}}\sum_{s \in \mathcal{S}}\Phi_{s}(t),
\end{multline}
where $\Phi_{i}$, $\Phi_{w}$, $\Phi_{k}$, and
$\Phi_{s}$ are penalty functions.

The conventional generation adjustment penalty function is
\begin{multline}
\label{eq:obj_conv_stage2}
\Phi_{i}(t) = c_{i}^{\mathrm{nl}} v_{i}^{}(t)
    + c_{i}^{\mathrm{su}} y_{i}^{}(t) \\
    + \sum_{b\in\mathcal{B}}\rho_{ib}
    \left\lbrack\delta\plus{p}_{ib}(t) + \delta\minus{p}_{ib}(t)\right\rbrack,
\end{multline}
where $\rho_{ib}$ is an incremental penalty on adjusting
generation dispatch.
The no-load and startup costs are included to ensure
changes in commitment are reflected in the objective.

The penalties arising from renewable generation curtailment
adjustments are given by
\begin{align}
    \label{eq:obj_wind_stage2}
    \Phi_{w}(t) &= \rho_{w} \left\lbrack \delta\plus{x}_{w}(t)
        + \delta\minus{x}_{w}(t) \right\rbrack \\
    \label{eq:obj_solar_stage2}
    \Phi_{k}(t) &= \rho_{k} \left\lbrack \delta\plus{x}_{k}(t)
        + \delta\minus{x}_{k}(t) \right\rbrack,
\end{align}
where $\rho_{w}$ and $\rho_{k}$ are incremental penalties on
adjusting wind and solar curtailment, respectively.

Lastly, the penalty associated with energy storage charging and discharging
is given by
%
%
%
\begin{equation}
    \label{eq:obj_ess_stage2}
    \Phi_{s}(t) = \rho_{s}\left\lbrack \delta{p}_{s}^{c}(t) + \delta{p}_{s}^{d}(t) \right\rbrack,
\end{equation}
where $\rho_{s}$ is an incremental penalty on adjusting ESS
injections.  This function is equivalent to penalizing the total
charging and discharging amounts because it differs only by a
constant.
%
%
For the complete energy storage model, refer to Section~\ref{sec:ess_model}.

\subsubsection{Generator output adjustment constraints}
The binary variable generation constraints and minimum up and
down time constraints in Stage~\num{2} are identical to those in Stage~\num{1};
however, the generator output constraints require modification.
Let $p_{ib}$ and $p_{i}$ be the block and unit outputs determined in
Stage~\num{1}.
%
%
Conventional generation output is then bounded by
\begin{align}
    \label{eq:st2_gen_con_1}
    \uvar{p}_{i}v_{i}(t) \leq
    p_{i}(t) &+ \delta \plus{p}_{i}(t) - \delta \minus{p}_{i}(t) \leq
    \ovar{p}_{i}v_{i}(t) \\
    \label{eq:st2_gen_con_2}
    0 \leq p_{ib}(t) &+ \delta \plus{p}_{ib}(t) - \delta \minus{p}_{ib}(t) \leq
    \ovar{p}_{ib}v_{i}(t),
\end{align}
for all ${b \in \mathcal{B}}$, ${i \in \mathcal{I}}$, and
${t \in \mathcal{T}}$.
As in \cref{eq:gen_con_1}, unit and block output are related by
\begin{multline}
    \label{eq:st2_gen_con_3}
    p_{i}(t) + \delta \plus{p}_{i}(t) - \delta \minus{p}_{i}(t) = \\
    \sum_{b\in \mathcal{B}}
        \left\lbrack
        p_{ib}(t) + \delta \plus{p}_{ib}(t) - \delta \minus{p}_{ib}(t)
        \right\rbrack,
\end{multline}
for all ${i \in \mathcal{I}}$ and ${t \in \mathcal{T}}$.

\subsubsection{Ramping constraints}

In Stage~\num{2}, conventional generation ramp rates are limited
such that
\begin{align}
    \nonumber
    -\uvar{R}_{i} &\leq p_{i}(t) + \delta\plus{p}_{i}(t)
    -\delta\minus{p}_{i}(t) \\
    &\quad- p_{i}(t-1) - \delta\plus{p}_{i}(t-1)
    + \delta \minus{p}_{i}(t-1) \leq \ovar{R}_{i},
\end{align}
for all ${i \in \mathcal{I}}$ and ${t \in \mathcal{T}}$.
This constraint follows from \cref{st1_ramp_constraints}.

\subsubsection{Nodal power balance and transmission constraints}
The nodal power balance constraints have the same structure
as~\cref{eq:nodal_power_bal}, with the exception
that the power terms are augmented to account for curtailments
and dispatch adjustments.
Since the objective of Stage~\num{2} is to determine
a minimal set of corrective actions required to
alleviate congestion,
we introduce transmission constraints
of the form
\begin{align}
    \label{eq:flow_limits}
    \uvar{F}_{l} \leq F_{l}(t) \le \ovar{F}_{l} \text{ for all } l \in \mathcal{L}_{m}
    \text{ and } t \in \mathcal{T},
\end{align}
where $\mathcal{L}_{m}$ denotes the set of monitored lines such that
${\mathcal{L}_{m} \subseteq \mathcal{L}}$.
For symmetric bidirectional flow limits, we have
${\uvar{F}_{l} = -\ovar{F}_{l}}$.

\subsubsection{Renewable generation curtailment adjustment constraints}
As in \cref{eq:st1_wind_curtail_con} and \cref{eq:st1_solar_curtail_con},
curtailment levels are bounded such that
\begin{align}
    0 &\leq x_{w}(t) + \delta\plus{x}_{w}(t) - \delta\minus{x}_{w}(t) \leq p_{w}(t) \\
    0 &\leq x_{k}(t) + \delta\plus{x}_{k}(t) - \delta\minus{x}_{k}(t) \leq p_{k}(t),
\end{align}
for all ${w \in \mathcal{W}}$, ${k \in \mathcal{K}}$, and ${t \in \mathcal{T}}$.

\subsubsection{Energy storage constraints}
\label{sec:ess_model}
The total energy storage charging and discharging amounts are given by
\begin{align}
    \label{eq:total_charge_def}
    p_{s}^{c}(t) &= p_{s}^{c_{0}}\<(t)
        + \delta\plus[c]{p}_{s}\<(t) - \delta\minus[c]{p}_{s}\<(t) \\
    \label{eq:total_discharge_def}
    p_{s}^{d}(t) &= p_{s}^{d_{0}}\<(t)
        + \delta\plus[d]{p}_{s}\<(t) - \delta\minus[d]{p}_{s}\<(t),
\end{align}
where $p_{s}^{c_{0}}$ is the initial charging schedule,
and $p_{s}^{d_{0}}$ the initial discharging schedule.
In Stage~\num{2} of the day-ahead formulation
${p_{s}^{c_0}(t) = p_{s}^{d_0}(t) = 0}$ for all ${t \in \mathcal{T}}$
because neither party has proposed a nonzero ESS injection schedule.
For a complete breakdown of these initial conditions by stage,
refer to \cref{tab:initial_ess_schedules} and
Section~\ref{sec:initial_ess_schedules}.

The charging and discharging decision variables are nonnegative
and bounded above such that
\begin{align}
    \label{eq:total_charge_bound}
    0 &\leq p_{s}^{c}(t) \leq \ovar{p}_{s}^{c}\nu_{s}(t) \\
    \label{eq:total_discharge_bound}
    0 &\leq p_{s}^{d}(t) \leq \ovar{p}_{s}^{d}\lbrack 1-\nu_{s}(t) \rbrack,
\end{align}
where $\ovar{p}_{s}^{c}$ is the maximum charging power of
the $s$th ESS, and $\ovar{p}_{s}^{d}$ the
maximum discharging power.
If the $s$th ESS is charging at time $t$, ${\nu_{s}(t)=1}$;
otherwise, ${\nu_{s}(t) = 0}$.
Constraints \cref{eq:total_charge_bound}
and \cref{eq:total_discharge_bound} prohibit simultaneous
charging and discharging.
Similarly, the charging adjustments
are nonnegative and bounded above such that
\begin{align}
    \label{eq:pos_charge_adjust}
    0 &\leq \delta\plus[c]{p}_{s}\<(t) \leq \lbrack\ovar{p}_{s}^{c}
        - p_{s}^{c_{0}}\<(t)\rbrack\psi_{s}(t) \\
    \label{eq:neg_charge_adjust}
    0 &\leq \delta\minus[c]{p}_{s}\<(t) \leq p_{s}^{c_0}\<(t)\lbrack 1-\psi_{s}(t) \rbrack,
\end{align}
where ${\psi_{s}(t) = 1}$ when the charging adjustment
of the $s$th ESS is positive, and ${\psi_{s}(t) = 0}$ otherwise.
For the discharging adjustments, we have
\begin{align}
    \label{eq:pos_discharge_adjust}
    0 &\leq \delta\plus[d]{p}_{s}\<(t) \leq \lbrack\ovar{p}_{s}^{d}
        - p_{s}^{d_{0}}\<(t)\rbrack\zeta_{s}(t) \\
    \label{eq:neg_discharge_adjust}
    0 &\leq \delta\minus[d]{p}_{s}\<(t) \leq p_{s}^{d_0}\<(t)\lbrack 1-\zeta_{s}(t) \rbrack,
\end{align}
where ${\zeta_{s}(t) = 1}$ when the discharging adjustment
of the $s$th ESS is positive, and ${\zeta_{s}(t) = 0}$ otherwise.
Constraints \cref{eq:total_charge_def}--\cref{eq:neg_discharge_adjust}
are enforced for all ${s \in \mathcal{S}}$ and ${t \in \mathcal{T}}$.

The formulation also includes a set of constraints that enable TEPO
to respect ESS state of charge limitations.
The difference equation that describes the SOC trajectory
is given by
\begin{align}
    \label{eq:ess_soc_def}
    E_{s}(t) &= E_{s}(t-1) + (\Delta\eta_{s}^{c})p_{s}^{c}(t) - (\Delta/\eta_{s}^{d})p_{s}^{d}(t),
\end{align}
for all ${s \in \mathcal{S}}$ and ${t \in \mathcal{T}}$,
where $\Delta$ is the step size parameter.
The charging efficiency of the $s$th ESS is denoted by
$\eta_{s}^{c}$, and the discharging efficiency by $\eta_{s}^{d}$.
The SOC is then bounded as follows:
\begin{align}
    \label{eq:ess_soc_bounds}
    \uvar{E}_{s}(t) &\leq E_{s}(t) \leq \ovar{E}_{s}(t),
\end{align}
for all ${s \in \mathcal{S}}$ and ${t \in \mathcal{T}}$,
where $\uvar{E}_{s}$ is the minimum SOC, and
$\ovar{E}_{s}$ the maximum SOC.
Additionally, we impose an equality constraint at the end of each
day such that
\begin{align}
    \label{eq:final_soc_target}
    E_{s}(t_{\<f}) &= E_{s}',
\end{align}
for all ${s \in \mathcal{S}}$.
The final time index of the day is $t_{\<f}$,
and the target SOC at ${t = t_{\<f}}$ is $E_{s}'$.
This approach brings each ESS to a predictable SOC at
the end of each day.

\subsubsection{Energy storage scheduling initial conditions by stage}
\label{sec:initial_ess_schedules}

The initial ESS charging and discharging schedules,
denoted respectively by $p_{s}^{c_0}$ and $p_{s}^{d_0}$,
vary with the optimization stage.
In Stage~\num{2} of the day-ahead formulation,
the initial schedules are zero-valued because neither party
has proposed a set of injections.
The schedules agreed upon at the conclusion of the day-ahead
framework, denoted by $\tilde{p}_{s}^{c}$ and $\tilde{p}_{s}^{d}$,
serve as the initial schedules for Stage~\num{2} of the hour-ahead
framework.
The schedule adjustments proposed by DEPO prior to Stage~\num{3}
are denoted by $\delta\tilde{p}_{s}^{c}$ and $\delta\tilde{p}_{s}^{d}$.
\Cref{tab:initial_ess_schedules} provides a complete breakdown
of the ESS scheduling initial conditions by stage.

\begin{table}[!t]
    \renewcommand{\arraystretch}{1.45}
    \caption{Initial charging and discharging conditions by stage}
    \label{tab:initial_ess_schedules}
    \centering
    \begin{tabular}{l l l l}
        \toprule
        Time frame & Stage & $p_{s}^{c_0}\<(t)$ & $p_{s}^{d_0}\<(t)$ \\
        \midrule
        Day-ahead & \textsf{\scriptsize CR 1} & 0, for all ${t \in \mathcal{T}}$ & 0, for all ${t \in \mathcal{T}}$ \\
        Day-ahead & \textsf{\scriptsize CR 2} & $\delta\tilde{p}_{s}^{c}\<(t)$ & $\delta\tilde{p}_{s}^{d}\<(t)$ \\
        Hour-ahead & \textsf{\scriptsize CR 1} & $\tilde{p}_{s}^{c}\<(t)$ & $\tilde{p}_{s}^{d}\<(t)$ \\
        Hour-ahead & \textsf{\scriptsize CR 2}
            & ${\tilde{p}_{s}^{c}\<(t) + \delta\tilde{p}_{s}^{c}(t)}$
            & ${\tilde{p}_{s}^{d}\<(t) + \delta\tilde{p}_{s}^{d}\<(t)}$ \\[0.3ex]
        \bottomrule
    \end{tabular}
\end{table}

\subsection{Stage~\num{3}, Coordinated congestion relief}
Stage~\num{3} is called the \textit{Coordinated congestion relief}
problem because it considers DEPO's proposed ESS schedule
(or hour-ahead adjustments).  The formulation is nearly the same as
in Stage~\num{2}, except the storage penalty function is given by
\begin{multline}
    \label{eq:obj_ess_stage3}
    \Phi_{s}(t) = \chi_{s}^{c_{0}}(t)\left\lbrack
                  \plus[c]\rho_{s}\delta\plus[c]{p}_{s}(t)
                + \minus[c]\rho_{s}\delta\minus[c]{p}_{s}(t)
                + \rho_{s}^{d_{\ast}}\delta\plus[d]{p}_{s}(t)\right\rbrack\\[0.5ex]
                + \chi_{s}^{d_{0}}(t)\left\lbrack
                  \plus[d]\rho_{s}\delta\plus[d]{p}_{s}(t)
                + \minus[d]\rho_{s}\delta\minus[d]{p}_{s}(t)
                + \rho_{s}^{c_{\ast}}\delta\plus[c]{p}_{s}(t)\right\rbrack,
\end{multline}
where $\chi_{s}^{c_{0}}$ and $\chi_{s}^{d_{0}}$ are indicator functions.
For charging, we have
\begin{align}
    \chi_{s}^{c_{0}}(t) &= \left\{\mkern-7mu
    \begin{array}{ll}
        1, & \text{if } {p}_{s}^{c_{0}}(t) > 0, \\
        0, & \text{otherwise,}
   \end{array}\right.
%
\end{align}
and $\chi_{s}^{d_{0}}$ is defined analogously for discharging.
Let $\plus[c]\rho_{s}$, $\minus[c]\rho_{s}$, $\plus[d]\rho_{s}$,
and $\minus[d]\rho_{s}$ be incremental penalties set by DEPO.
These penalties reflect the value that DEPO places on maintaining a
particular ESS injection schedule.  The parameters
$\rho_{s}^{c_{\ast}}$ and $\rho_{s}^{d_{\ast}}$ are similar, but they
allow DEPO to ascribe different incremental penalties when TEPO
reverses the ESS charging action.  Here these incremental penalties
are in units of \num{1}/\si{\MWh}, although in theory they could be
viewed as prices.
%
%
The output of Stage~\num{3} fixes the ESS injections and, by extension,
the net load at the storage buses.
After completing this stage, TEPO sends DEPO
the mitigation needs report.
Based on this information, DEPO computes and returns the final
combined ESS schedule.

\subsection{Stage~\num{4}, Post-mitigation unit commitment}
Stage~\num{4} is called the \textit{Post-mitigation unit commitment}
problem.  In it, TEPO solves a network-constrained unit commitment and
economic dispatch.  The constraints are the same as in Stage~\num{1},
except the line flow limits in \cref{eq:flow_limits} are also
enforced.  The objective function is identical to
\cref{eq:stage1_objective}.  This problem treats the ESS injection
schedules determined in Stage~\num{3} as inputs and solves for the
least-cost generation schedule.  The output of Stage~\num{4} is the
optimal commitment and dispatch considering ESS injections and
transmission capacity constraints.

\section{3-bus case study}
\label{sec:3bus_case}
\begin{figure}[!t]
    \centering
    \includegraphics[width=3.2in]{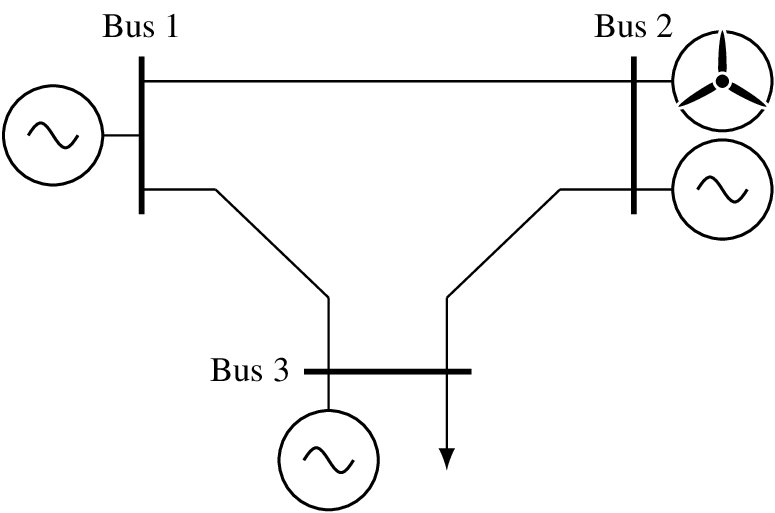}
    \caption{3-bus test system.}
    \label{fig:3bus_system}
\end{figure}
\begin{table}[!t]
    \renewcommand{\arraystretch}{1.15}
    \caption{3-bus system line data}
    \label{tab:3bus_lines}
    \centering
    \begin{tabular}{l r r}
        \toprule
        Line & Reactance (\si{\pu}) & Max.\ flow (\si{\MW}) \\
        \midrule
        1--2 & 0.13 & 50 \\
        1--3 & 0.13 & 50 \\
        2--3 & 0.13 & 25 \\
        \bottomrule
    \end{tabular}
    \vspace*{1.25em}
    \caption{3-bus system generator data}
    \label{tab:3bus_gens}
    \centering
    \begin{tabular}{l r r r r r}
        \toprule
        Unit & Bus & Min.\ gen.\ & Max.\ gen.\ & Start-up cost & Inc.\ cost \\
             & & (\si{\MW})   & (\si{\MW})   & (\$) & (\$/\si{\MWh}) \\
        \midrule
        1 & 1 & 10 & 100 & 100 & 30 \\
        2 & 2 & 10 & 100 & 100 & 40 \\
        3 & 3 & 10 &  50 & 100 & 20 \\
        \bottomrule
    \end{tabular}
    \vspace*{1.25em}
    \caption{3-bus results summary}
    \label{tab:3bus_results}
    \centering
    \begin{tabular}{l r r r r}
        \toprule
        Location & Total cost & Gen.\ cost & Spill.\ cost & Wind spill.\ \\
        & (\$) & (\$) & (\$) & (\si{\MWh}) \\
        \midrule
        No ESS & \num{25090} & \num{24470} & \num{620} & \num{31}   \\
        Bus 1  & \num{25047} & \num{24427} & \num{620} & \num{31}   \\
        Bus 2  & \num{24692} & \num{24272} & \num{420} & \num{21}   \\
        Bus 3  & \num{24642} & \num{24184} & \num{458} & \num{22.9} \\
        \bottomrule
    \end{tabular}
\end{table}
To illustrate the formulation presented in \cref{sec:formulation}, a
case study based on a small test system was
developed. \Cref{fig:3bus_system} shows the system of interest.  The
transmission line parameters are given in \cref{tab:3bus_lines}, and
the generator data in \cref{tab:3bus_gens}.  The system load is
concentrated at Bus~3 and reaches a peak of \SI{110}{\MW}.
Here we explore the effect of siting a
\SI{5}{\MW}/\SI{10}{\MWh} battery at one of the buses in the
system.  The full ESS capacity was made available to TEPO, reflecting
no local service provision.  For simplicity, the constraints on
generator ramp rates and minimum up and down times were relaxed, and
all of the conventional generators were initially scheduled off.  The
large-scale case study in \cref{sec:case_study_and_demo} considers all
of the constraints.

\begin{figure}[!t]
    \centering
    \includegraphics[width=3.4in]{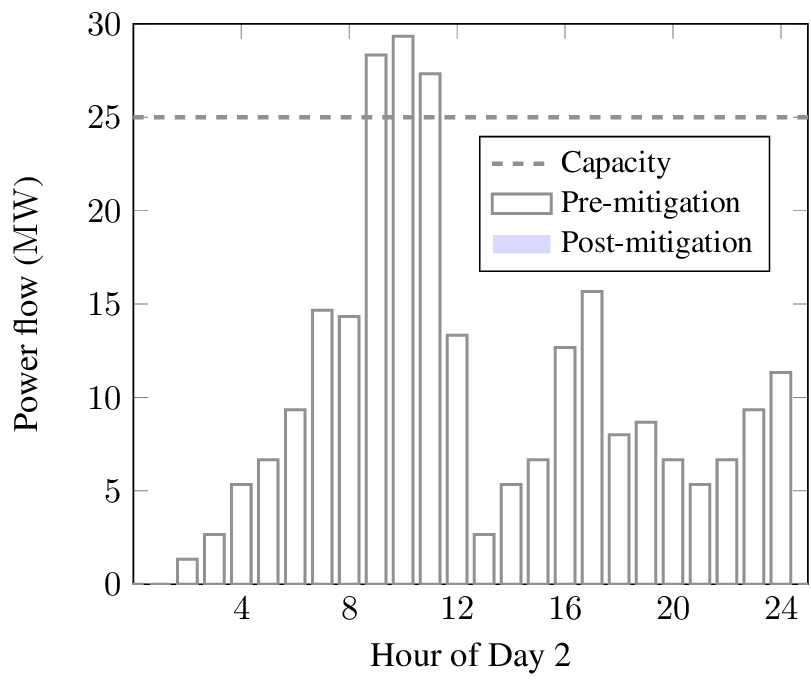}
    \caption{3-bus case study: Congestion mitigation on Line 2--3.}
    \label{fig:3bus_congestion}
\end{figure}

Over the operating day of interest, the transmission line connecting
the wind generation to the load center becomes congested.
Because Unit~3 is located next to the load and has the lowest
incremental cost, its output is maximized during the period of
congestion.  Unit~2 is not committed.  Hence, in the absence of energy
storage, wind curtailment is the primary mechanism used to reduce
congestion on Line~\mbox{2--3}.  The amount of wind energy curtailed
can be reduced by installing energy storage, depending on where the
ESS is sited.  \Cref{fig:3bus_congestion} shows the real power flow on
Line~\mbox{2--3} in the pre- and post-mitigation stages when the ESS
is sited at Bus~3.  As shown in \cref{tab:3bus_results}, placing the
ESS here yields the lowest overall operation cost and reduces the wind
curtailment by roughly \SI{26}{\percent} over the case with no
storage. Intuitively, the ESS location that yields the lowest
curtailment is Bus 2, next to the wind farm. In this particular case,
the operation cost is slightly higher than when the ESS is sited near
the load due to differences in generation dispatch.  This example
demonstrates that the scheduling methodology promotes congestion
relief, cost savings, and improved renewable energy integration.

\section{Large-scale case study and demonstration}
\label{sec:case_study_and_demo}

\begin{figure}[!t]
    \centering
    \includegraphics[width=3.4in]{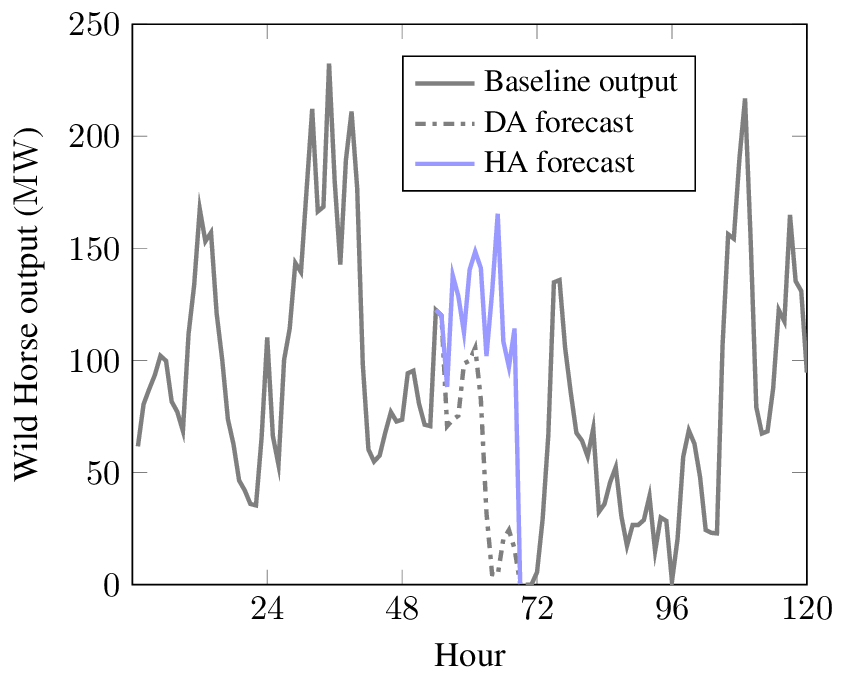}
    \caption{Wild Horse wind ramp.}
    \label{fig:wild_horse_ramp}
\end{figure}

To demonstrate the scalability of the framework, we developed a
large-scale case study based on a realistic representation of the
Pacific Northwest region of the United States.  The system model,
summarized in \cref{tab:system_description}, was built from a subset
of the WECC 2024 Common Case~\cite{WECC:14}.  This production cost
modeling data set projects how the generation mix of the Western
Interconnection will change over time.  The focal point of the case
study was an actual \SI{2}{\MW}/\SI{1}{\MWh} battery energy storage
system located at SnoPUD's Hardeson substation in Everett, WA.  By
itself, this amount of capacity is insufficient to significantly
reduce transmission congestion; therefore, it was represented in the
model as a \SI{200}{\MW}/\SI{100}{\MWh} ESS.  This served to better
illustrate the capabilities of the coordination framework given
adequate resources.

\begin{figure}[!t]
    \centering
    \hspace*{-2.5mm}
    \includegraphics[width=3.375in]{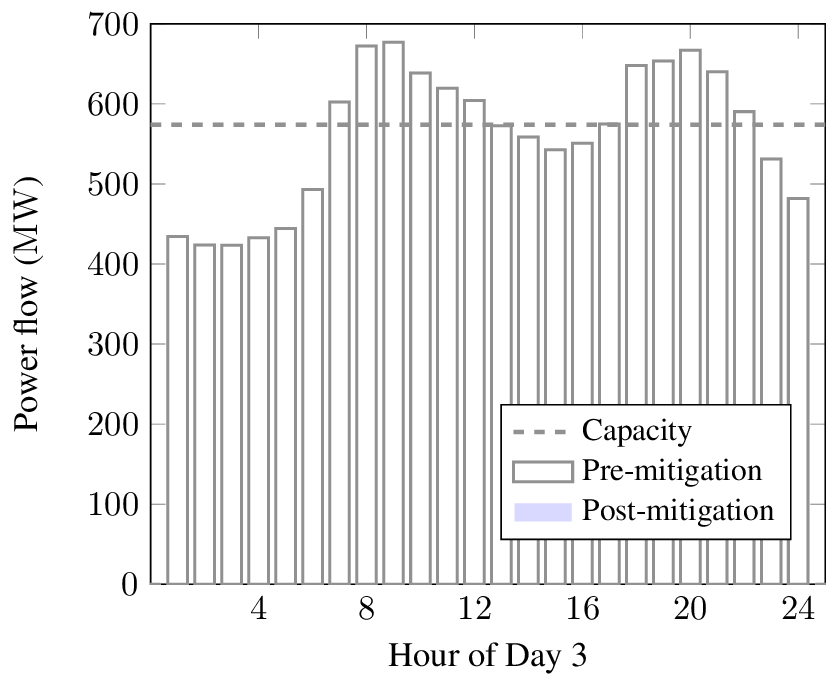}
    \caption{Congestion mitigation on Line 1716, Blue Lake--Troutdale 230 kV.}
    \label{fig:congestion_mitigation}
\end{figure}

This case study examines a scenario in which there are substantial
changes in the wind forecast between the day-ahead and hour-ahead
frameworks.  Specifically, we consider a multi-hour wind ramp at the
Wild Horse wind farm (\SI{273}{\MW} capacity) near Ellensburg,
WA.  The day-ahead and hour-ahead wind power forecasts at Wild Horse
are shown in \cref{fig:wild_horse_ramp}.
%
%
During this operating day, the \num{1.5} mile Blue Lake--Troutdale
\SI{230}{\kV} transmission line near Portland, OR is congested for
roughly \num{11} hours.
That is, the solution to a standard unit commitment and economic
dispatch causes violations of the short term line rating.
\cref{fig:congestion_mitigation} shows the real power flow on this
line in both the pre- and post-mitigation stages.
%

\begin{figure}[!t]
    \centering
    \hspace*{2.0mm}
    \includegraphics[width=3.4in]{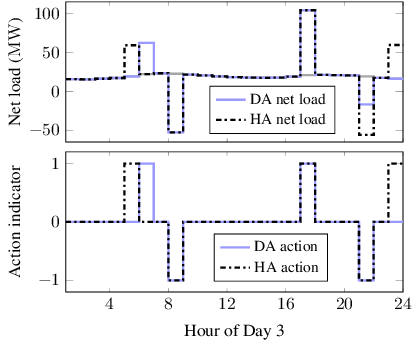}
    \caption{Net load at the SnoPUD energy storage bus.}
    \label{fig:ess_net_load}
\end{figure}

The optimal net injections for the ESS are shown in
\cref{fig:ess_net_load}.
An action indicator value of \num{1} implies the ESS is charging, and
a value of \num{-1} implies the ESS is discharging.  The injections
shown in \cref{fig:ess_net_load} correspond to the two distinct blocks
of time when the Blue Lake--Troutdale line is congested.  Immediately
prior to each congested period the battery is charged so that it may
discharge at the correct moment to mitigate congestion.  In
conjunction with modestly redispatching some conventional generation
units, this strategy is effective in satisfying the static system
security constraints.  A notable finding of this study is that the
storage system was able to help mitigate congestion on a transmission
line despite being roughly \num{170} miles away.

In \cref{fig:ess_net_load}, the ESS injections determined in the
day-ahead and hour-ahead frameworks do not fully overlap.  In this
case, the mismatch between the two is expected because of the
substantial change in the system operating condition.  Effectively,
the solution provided by the day-ahead framework is no longer optimal
because of the changes in the wind forecast.  The hour-ahead framework
allows TEPO to proactively plan for the wind ramp.  The charging
activity at the end of the day brings the ESS back to its target SOC,
$E_{s}'$.

\begin{table}[!t]
    \renewcommand{\arraystretch}{1.15}
    \caption{Large-scale test system information}
    \label{tab:system_description}
    \centering
    \begin{tabular}{l r}
        \toprule
        Component & Quantity \\
        \midrule
        Buses & \num{2764} \\
        Branches & \num{3318} \\
        Fixed generators & \num{440} \\
        Controllable generators & \num{38} \\
        Wind farms & \num{73} \\
        Solar plants & \num{5} \\
        Energy storage systems & \num{1} \\
        \bottomrule
    \end{tabular}
\end{table}

\subsection{Demonstration}
\label{sec:demo}

The framework presented in Sections \ref{sec:methodology} and
\ref{sec:formulation} was implemented within SnoPUD's production power
scheduling environment.  TEPO was implemented in GAMS using the IBM
CPLEX solver~\cite{GAMS:17, San:04}.  DEPO implements an
\mbox{OpenADR} Virtual End Node that communicates with TEPO via a
cloud-hosted XMPP server~\cite{McP:11, Hor:10}.  Upon completion of
the optimization procedure, DEPO supplies ESS scheduling
recommendations to a human operator.  The large-scale case study
outlined above served as the basis for the live demonstration.  To
highlight the transmission-side impact, the ESS was not called upon
for local service provision during the demonstration.  Over the
operating day with the wind ramp, the maximum absolute difference
between the simulated and actual schedules was \SI{70}{\kW} for the
\SI{2}{\MW} battery system.  The scheduling differences largely
reflect approximation errors in the SOC tracking
constraints \cref{eq:ess_soc_def}--\cref{eq:final_soc_target}.

\section{Benchmark with one-shot optimization}
\label{sec:benchmark}
Here we explore the possibility that, given sufficient centralized
coordination, TEPO could solve a master problem encompassing all
constraints and distribution system requirements.  This hypothetical
formulation would take the form of a one-shot optimization rather than
a sequence of linked stages.  The design of this master problem would
need to carefully specify what information needs to be exchanged and
when between TEPO and the DEPO instances.  This one-shot formulation
could potentially be decomposed preserving some degree of independence
and information privacy by allowing each DEPO instance to solve its
own subproblem.  The complete development of a master problem that
meets the above criteria is outside the scope of this paper; however,
when the DEPOs make their storage capacity fully available to TEPO, as
in the case study from \cref{sec:case_study_and_demo}, the master
problem effectively becomes a network-constrained unit commitment
(NCUC) with energy storage constraints.  Thus, we present a benchmark
comparison of the operating costs incurred in Stage~\num{4} of the
multistage framework versus the augmented NCUC, henceforth denoted as
NCUC+.  The purpose of this comparison is to provide a rough estimate
of the cost of not solving the problem in a fully centralized manner.

\subsection{Augmented network-constrained unit commitment}
The tuple of decision variables in NCUC+ is
\begin{equation}
    \nonumber
    \xi : \left(v_{i}, y_{i}, z_{i}, \nu_{s}, p_{ib},
        p_{i}, p_{s}^{c}, p_{s}^{d}, x_{k}, x_{w}, \theta_{n} \right),
\end{equation}
where $p_{s}^{c}$ and $p_{s}^{d}$ are the charging and discharging
schedules of the $s$th ESS, respectively.  As in
\cref{eq:total_charge_bound}, $\nu_{s}$ is a binary variable that
prevents simultaneous charging and discharging.
The remaining decision variables are defined as in Stage~\num{1}.

The objective function of NCUC+ is:
\begin{multline}
    \label{eq:stage5_objective}
    f(\xi) = \sum_{t \in \mathcal{T}}\sum_{i \in \mathcal{I}}C_{i}(t)
    + \sum_{t \in \mathcal{T}}\sum_{w \in \mathcal{W}}C_{w}(t) \\
    + \sum_{t \in \mathcal{T}}\sum_{k \in \mathcal{K}}C_{k}(t)
    + \sum_{t \in \mathcal{T}}\sum_{s \in \mathcal{S}}C_{s}(t),
\end{multline}
where $C_{i}$, $C_{w}$, and $C_{k}$ are defined as in
\cref{eq:stage1_objective}.  The energy storage utilization costs
are
\begin{align}
    C_{s}(t) = m_{s}\left\lbrack{p_{s}^{c}(t) + p_{s}^{d}(t)}\right\rbrack,
\end{align}
where charging and discharging are priced symmetrically for
simplicity.  In order to make a fair comparison, the incremental cost of
storage utilization was set to match the incremental penalty on ESS
injections from \cref{eq:obj_ess_stage2}, \ie, ${m_{s} = \rho_{s}}$.
The minimum operating cost objective in \cref{eq:stage5_objective}
results in the available storage capacity being used for a variety of
transmission services, such as temporal arbitrage, rather than purely
for congestion relief.  The constraints of NCUC+ include those from
Stage~\num{4} and a simplified energy storage model based on
(\ref{eq:total_charge_def})--(\ref{eq:final_soc_target}).

\subsection{Quantitative comparison}
\begin{table}[!t]
    \renewcommand{\arraystretch}{1.15}
    \caption{Benchmark results summary}
    \label{tab:benchmark_summary}
    \centering
    \begin{tabular}{l r r r r}
        \toprule
        Problem & Storage cap. & Gen. cost & Cost reduction & Reduct. gap \\
        & (\si{\MWh}) & (k\$) & (\$) & (\si{\percent}) \\
        \midrule
        NCUC+ & \num{0} & \num{1740.93} & -- & -- \\
        NCUC+ & \num{100} & \num{1735.96} & \num{4964.27} & \num{0.0} \\
        Stage~\num{4} & \num{100} & \num{1735.97} & \num{4955.09} & \num{0.2} \\
        \bottomrule
    \end{tabular}
\end{table}

The large-scale test system described in
\cref{sec:case_study_and_demo} was used for the benchmark comparison.
As in the case study, the focal point of this analysis was operating
day~\num{3}.  Each data point in the comparison corresponds to a
scheduling method paired with a given amount of storage capacity.  To
determine a baseline, we ran NCUC+ with no energy storage.  When there
is no storage in the system, NCUC+ is mathematically equivalent to
Stage~\num{4}, \ie, they return exactly the same solution.  Then the
\SI{200}{\MW}/\SI{100}{\MWh} battery was reinserted and the multistage
and one-shot frameworks were compared.
For all optimization runs, a relative MILP gap of \SI{0.1}{\percent}
was used as the stopping criterion, \ie, \texttt{optcr=0.001} in GAMS.
The results are summarized in \cref{tab:benchmark_summary}.

For the case with no energy storage, the generation costs are
\SI{1740926}[\$]{} over the day of interest.  With the battery in the
system, the generation costs incurred by NCUC+ are reduced by
\SI{4964.27}[\$]{}.  This difference is attributable to the ESS, which
NCUC+ uses to perform a variety of transmission-side services.  In
contrast, the generation costs incurred in Stage~\num{4} of the
multistage framework decrease by \SI{4955.09}[\$]{} in relation to the
case with no storage.  Although the reduction in operating cost is
slightly smaller than with NCUC+, the difference is only about
\SI{0.2}{\percent}.  This difference, indicated in
\cref{tab:benchmark_summary} as the cost \textit{reduction gap}, can
be partially attributed to the fact that the one-shot formulation uses
the ESS to perform multiple transmission-side services while the
multistage framework does not.  The cost reduction gap is
sensitive to the incremental penalties employed in the congestion
relief stages.  For instance, if the incremental penalty on wind
curtailment is set much lower than the actual incremental cost, \ie,
${\rho_{w} \ll m_{w}}$, the congestion relief stages may produce storage
injection schedules that presume too much wind curtailment and
therefore have a larger gap.  In this analysis, the incremental
penalties were set to reflect the actual costs, \eg,
${\rho_{w} = m_{w}}$.

From an economic perspective, the cost reduction gap is indicative of
how the social welfare declines when TEPO and DEPO act in their own
self-interest.
This benchmark comparison indicates that there is indeed a price to be
paid for allowing TEPO and DEPO to act independently, but that price
appears to be small when the incremental penalties in Stages~\num{2}
and \num{3} are set appropriately.  Eliminating this effect entirely
would require full centralized coordination and/or a carefully crafted
decomposition of a suitable master problem.

\section{Conclusion}
\label{sec:conclusion}

This paper addresses the problem of how to share energy storage
capacity among transmission and distribution entities.  It describes
and demonstrates a method for coordinating transmission-level
congestion relief with local, distribution-level services in systems
that lack centralized markets.  A weighted ${\ell_{1}\text{-norm}}$
objective determines a minimal set of corrective actions required to
alleviate transmission congestion.  This work could be readily
extended to accommodate other system-wide objectives, such as
frequency regulation.  Future work will explore the effect of line
losses when sharing energy storage capacity among transmission and
distribution entities.  Finally, another interesting avenue of
research will be the use of mathematical decomposition techniques
(\eg, Dantzig-Wolfe decomposition) to analyze the coordination of
energy storage services between transmission and distribution within
a centralized environment.

\ifCLASSOPTIONcaptionsoff
  \newpage
\fi



\bibliographystyle{IEEEtran}
\bibliography{IEEEabrv,./bib/uw1e_journal}


%
%
%
%
\vskip 0.075\baselineskip plus -1fil
\begin{IEEEbiographynophoto}{Ryan Elliott}
is a Ph.D.\ candidate in the Department of Electrical Engineering at
the University of Washington.  His research focuses on renewable
energy integration, wide-area measurement systems, and power system
operation and control. Prior to pursuing a Ph.D., he was with the
Electric Power Systems Research Department at Sandia National
Laboratories from 2012--2015. While at Sandia, he served on the WECC
Renewable Energy Modeling Task Force, leading the development of the
WECC model validation guideline for central-station PV plants.  In
2017, he earned an R\&D 100 Award for his contributions to the design
of a real-time damping control system using PMU feedback.  Ryan
received the M.S.E.E.\ degree from the University of Washington in
2012.
\end{IEEEbiographynophoto}
\vskip -1\baselineskip plus -1fil
\begin{IEEEbiographynophoto}{Ricardo Fern\'{a}ndez-Blanco}
is a Postdoctoral Researcher at the University of Malaga, Spain.
His research interests include the fields of operations and economics
of power systems, smart grids, bilevel programming, hydrothermal
coordination, electricity markets, and the water-energy nexus.
Previously, he was a Postdoctoral Researcher at the University of
Washington, Seattle, WA, USA, and later a Scientific/Technical Project
Officer in the Knowledge for the Energy Union Unit at the Joint
Research Centre (DG JRC) of the European Commission, Petten, The
Netherlands.  Ricardo received the Ingeniero Industrial degree and the
Ph.D.\ degree in electrical engineering from the Universidad de
Castilla-La Mancha, Ciudad Real, Spain, in 2009 and 2014,
respectively.
\end{IEEEbiographynophoto}
\vskip -1\baselineskip plus -1fil
\begin{IEEEbiographynophoto}{Kelly Kozdras}
is the Distributed Energy Resource Planner at Puget Sound Energy
(PSE).  Previously while at PSE she provided engineering services for
PSE-owned generating facilities and pilot projects in utility and
customer-scale battery storage.  She earned her M.S.E.E.\ degree from
University of Washington in 2016, completing her thesis on modeling
and analysis of a microgrid containing hydro and battery
storage. Prior to that, Kelly worked in various capacities as an
electrical engineer and task manager on infrastructure projects around
the USA as well as at South Pole Station, Antarctica.  She earned a
B.S.E.E.\ from Rose-Hulman Institute of Technology in 1999, and is a
licensed professional engineer in the state of Washington.
\end{IEEEbiographynophoto}
\newpage
\begin{IEEEbiographynophoto}{Josh Kaplan}
is an aspiring carpenter and avid outdoorsman based in Seattle, WA.
Formerly, he was a founding employee of 1Energy Systems, an energy
storage software company.  At 1Energy, which became Doosan GridTech
following acquisition, he was responsible for software design and
project management on several energy storage deployments and research
projects.  Prior to his work in the energy storage field, he worked as
a software engineer at Microsoft, primarily on developer tools
products. Josh holds a B.A.\ degree in Political Science from the
University of Washington.
\end{IEEEbiographynophoto}
\vskip 1\baselineskip plus -1fil
\begin{IEEEbiographynophoto}{Brian Lockyear}
is a Principal Software Development Engineer with Doosan GridTech. His
interests include green energy in both design and implementation and
in applications of artificial intelligence in energy control
systems. He holds a Ph.D.\ in computer science from the University of
Washington, a M.Arch.\ from the University of Oregon, and a
B.S.E.E.\ from Oregon State University. Prior to Doosan, he worked for
Synopsys, Tera Computer, and NASA.
\end{IEEEbiographynophoto}
\vskip 1.62\baselineskip plus -1fil
\begin{IEEEbiographynophoto}{Jason Zyskowski}
is the Senior Manager of Planning, Engineering and Technical Services
at the Snohomish County PUD.  He was the project manager for the
District's first energy storage system deployment and has contributed
to numerous renewable generation and automation upgrade projects.
Jason has been with the PUD since receiving the B.S.E.E.\ degree from
the University of Washington in 2004.  His professional experience
includes time in the Transmission, System Protection and Substation
Engineering departments.  He is a registered professional engineer in
the state of Washington.
\end{IEEEbiographynophoto}
\vskip 1.62\baselineskip plus -1fil
\begin{IEEEbiographynophoto}{Daniel Kirschen}
is the Donald W. and Ruth Mary Close Professor of Electrical
Engineering at the University of Washington. His research focuses on
the integration of renewable energy sources in the grid, power system
economics and power system resilience. Prior to joining the University
of Washington, he taught for 16 years at The University of Manchester
(UK). Before becoming an academic, he worked for Control Data and
Siemens on the development of application software for utility control
centers. He holds a Ph.D.\ from the University of Wisconsin-Madison and
an Electro-Mechanical Engineering degree from the Free University of
Brussels (Belgium). He is the author of two books.
\end{IEEEbiographynophoto}
\vfill

%
%




\end{document}